\documentclass[a4paper,12pt,reqno]{amsart}
\usepackage{latexsym}
\usepackage{amsmath,amsthm,amssymb,amscd}
\usepackage{fullpage}
\usepackage{xcolor}
\usepackage{csquotes}
\usepackage{parskip}
\usepackage{mathtools}
\mathtoolsset{showonlyrefs}

\usepackage[activate={true,nocompatibility},final,tracking=true,kerning=true,spacing=true]{microtype}
\usepackage[pdfencoding=unicode,pdftex,bookmarks=true,bookmarksnumbered, pagebackref]{hyperref}
\definecolor{citecolour}{rgb}{0.0, 0.0, 0.8}

\colorlet{linkcolour}{green!50!black}
\hypersetup{colorlinks,breaklinks,
		linkcolor=linkcolour,citecolor=citecolour,
		filecolor=linkcolour, urlcolor=linkcolour}
\usepackage{version}

\newtheorem{prevtheorem}{Theorem}

\newtheorem{prevcorollary}{Corollary}

\newtheorem{prevlemma}{Lemma}

\newtheorem{theorem}{Theorem}

\newtheorem*{theoremno}{Theorem}
\newtheorem{problem}[theorem]{Problem}
\newtheorem{question}[theorem]{Question}

\theoremstyle{remark}

\numberwithin{equation}{section}
\usepackage[shortlabels]{enumitem}
\begingroup
 \makeatletter
 \@for\theoremstyle:=definition,remark,plain\do{%
 \expandafter\g@addto@macro\csname th@\theoremstyle\endcsname{%
 \addtolength\thm@preskip\parskip
 }%
 }

\DeclareMathOperator{\X}{\mathfrak{X}}

\DeclareMathOperator{\N}{\mathbf{N}}

\DeclareMathOperator{\cent}{\mathbf{C}}

\DeclareMathOperator{\nsup}{nsup}

\DeclareMathOperator{\aut}{\mathbf{Aut}}
\DeclareMathOperator{\gl}{GL}
\DeclareMathOperator{\zent}{\mathbf{Z}}
\DeclareMathOperator{\fratt}{\mathbf{\Phi}}

\newcommand{\R}{\mathbb{\R}}

\includeversion{comment}
\renewcommand{\leq}{\leqslant}
\renewcommand{\geq}{\geqslant}
\newcommand{\md}[1]{\,\left(\textnormal{mod}\ #1\right)}

\newenvironment{proofofA}{{\bf {Proof of Theorem \ref{Thm:A}.} }}{\hfill $\blacksquare$ \\}
\newenvironment{proofofB}{{\bf {Proof of Theorem \ref{Thm:B}.} }}{\hfill $\blacksquare$ \\}
\newenvironment{proofofC}{{\bf {Proof of Theorem \ref{Thm:C}.} }}{\hfill $\blacksquare$ \\}
\newenvironment{proofof}{{\bf {Proof.} }}{\hfill $\blacksquare$ \\}

\begin{document}
\title{Subgroup congruences for groups of prime power order}
\author{Stefanos Aivazidis}
\address{Department of Mathematics and Applied Mathematics, University of Crete, Voutes Campus, 70013 Heraklion, Crete, Greece.}
\email{s.aivazidis@uoc.gr}
\author{Maria Loukaki}
\address{Department of Mathematics and Applied Mathematics, University of Crete, Voutes Campus, 70013 Heraklion, Crete, Greece.}
\email{mloukaki@uoc.gr}

\subjclass[2010]{20D15 (20D60)}
\begin{abstract}
Given a $p$-group $G$ and a subgroup-closed class $\X$, we associate 
with each $\X$-subgroup $H$ certain quantities which count $\X$-subgroups
containing $H$ subject to further properties. 
We show in Theorem~\ref{Thm:A} that each one of the said quantities is always $\equiv 1 \md{p}$
if and only if the same holds for the others.
In Theorem~\ref{Thm:B} we supplement the above result by focusing on normal $\X$-subgroups
and in Theorem~\ref{Thm:C} we obtain a sharpened
version of a celebrated theorem of Burnside relative
to the class of abelian groups of bounded exponent.
Various other corollaries are also presented.
\end{abstract}
\maketitle
\tableofcontents

\section{Introduction and main results}
The starting point of the work reported here was a theorem
announced by Lior Yanovski. Yanovski wrote an e-mail to the 
\enquote{group-pub-forum} online discussion group enquiring if anyone was
aware of an already existing proof for the following:

\begin{quote}
If $G$ is a finite $p$-group, then the number of maximal
abelian subgroups of $G$ is $\equiv 1 \md{p}$.
\end{quote}

Some discussion ensued, but ultimately no one claimed to have
seen this result before. Marty Isaacs asked Yanovski to share
his proof with the other pubbers and Yanovski obliged. (This result can now be found in \cite{Yanovski}.)
Later, Isaacs went on to generalise Yanovski's theorem
by replacing the condition \enquote{maximal abelian subgroups}
with the stronger 
\enquote{maximal abelian subgroups of exponent $\leq p^k$, where $p^k>2$}. (This generalisation is the main result in \cite{isaacs}.)

The present paper is the culmination of our desire to look deeper
both into Yanovski's method of proof by M\"{o}bius inversion 
and to provide a unified framework for counting 
subgroups with specified properties within $p$-groups. 
In particular, we work with general classes of groups $\X$ which are
subgroup-closed.
Our first main result is the following.

\begin{prevtheorem}\label{Thm:A}
Let $G$ be a finite $p$-group
and let $\X$ be a subgroup-closed class of finite groups.
Then the following are equivalent:
\begin{enumerate}[label={\upshape(\Alph*)}]
\item \label{conj:a}For each $\X$-subgroup $H$ of $G$, 
if $\mathfrak{a}_G(H) \neq 0$ then $\mathfrak{a}_G(H) \equiv 1 \md{p}$.
\item \label{conj:b}For each $\X$-subgroup $H$ of $G$ we have that 
$\mathfrak{b}_G(H) \equiv 1 \md{p}$.
\item \label{conj:c}For each $\X$-subgroup $H$ of $G$ and for all integers $k$ 
such that $|H| \leq p^k \leq \sup_G(H)$ 
we have that $\mathfrak{c}_{k, G}(H) \equiv 1 \md{p}$.
\end{enumerate}
\end{prevtheorem}

Here $\mathfrak{a}_G(H)$ is the number of $\X$-subgroups
which contain $H$ as a maximal subgroup and $\mathfrak{b}_G(H)$
is the number of $\X$-subgroups which are maximal subject to
containing $H$.
In general, there may exist several subgroups which are maximal
subject to containing $H$ and some of these certainly may have
different orders. We write $\sup_G(H)$ to denote the smallest
such order. Our notation $\mathfrak{c}_{k, G}(H)$ stands for the
number of $\X$-subgroups which contain $H$ and lie at a given level,
i.e. have fixed order $p^k$.

Our second main theorem acts as a supplement 
to the first and reads:

\begin{prevtheorem}\label{Thm:B}
Let $G$ be a finite $p$-group
and let $\X$ be a subgroup-closed class of finite groups.
The group $G$ satisfies one (and thus all) of the conditions \ref{conj:a}, \ref{conj:b}, \ref{conj:c} in Theorem~\ref{Thm:A} if and only if $G$ satisfies the following two conditions: 
\begin{enumerate}[label={\upshape(\arabic*)}]
\item\label{Thm2a} Every maximal normal $\X$-subgroup of $G$ is also a maximal $\X$-subgroup of $G$. 
\item\label{conj:d} For each normal $\X$-subgroup $H$ of $G$ and for all integers $k$ 
such that $|H| \leq p^k \leq \nsup_G(H)$ 
we have that $\mathfrak{c}_{k, G}(H) \equiv 1 \md{p}$.
\end{enumerate}
\end{prevtheorem}

Our notation $\nsup_G(H)$ is the appropriate analogue of $\sup_G(H)$
for normal subgroups.
All notational conventions and details of our set-up will be reiterated
and further explained in the next section.

In the last section of our paper we demonstrate how our main theorems may be
used jointly with certain auxiliary results to recover the previously mentioned results of Yanovski and Isaacs, as well as a classical
result due to Miller which guarantees the existence of abelian subgroups
in a certain range. 
Moreover, we amplify a celebrated theorem due to Burnside. In particular, our last main result, 
which may reasonably be said to have some independent interest,
asserts the following.

\begin{prevtheorem}\label{Thm:C}
Suppose that $G$ is a group of order $p^n$ and exponent $p^e$. 
Let $A \leq G$ be maximal among the normal abelian subgroups of $G$ 
having exponent $ p^s$, where $s \geq 2$. 
If $|A| = p^r$ then 
\[
n \leq \frac{1}{2}(r - s + 1)(r - s + 2e - 2) +r \, ,
\]
and thus 
\[
r \geq -(e-s +1/2) + \sqrt{(e-s+1/2)^2 +2n}
\]
\end{prevtheorem}

\section{Proofs of Main Results}
We denote by $\X$ a class of finite groups, i.e. a collection of groups
which is closed under isomorphisms. We assume, further, that $\X$
is closed under the taking of subgroups. We call $G$ an $\X$-group
if $G \in \X$ in which case all subgroups of $G$ are $\X$-groups
and thus $\X$-subgroups of $G$. We call a subgroup $H$ of $G$
a maximal $\X$-subgroup if there exists no $K \leq G$ such that
$K \in \X$ and $H < K$.

 Given a prime $p$, 
a finite $p$-group $G$
and an $\X$-subgroup $H$ of $G$, we define the following quantities
\begin{itemize}
\item $\mathfrak{a}_G(H) \coloneqq 
\left\lvert \{H < K \leq G : K \in \X \text{ and } |K:H|=p \} \right\rvert$
so that $\mathfrak{a}_G(H)$ is the number of minimal overgroups of $H$
which are $\X$-groups ;
\item $\mathfrak{b}_G(H) \coloneqq 
\left\lvert \{H \leq K \leq G : \text{$K$ is a maximal $\X$-subgroup} \} \right\rvert$
so that $\mathfrak{b}_G(H)$ is the number of maximal $\X$-subgroups
of $G$ which contain $H$;
\item $\mathfrak{c}_{k, G}(H) \coloneqq 
\left\lvert \{H \leq K \leq G : K \in \X \text{ and } |K|=p^k \} \right\rvert$ 
so that $\mathfrak{c}_{k, G}(H)$ is the number of 
$\X$-subgroups of $G$ which contain $H$ and have order $p^k$;
\item $\sup_G(H) \coloneqq \min\{|K| : H \leq K
\text{ and $K$ is a maximal $\X$-subgroup}\}$; 
in other words, $\sup_G(H)$ is the \enquote{lowest level}
at which a maximal $\X$-subgroup of $G$ containing $H$ occurs.
\item $\nsup_G(H) \coloneqq \min\{|K| : H \leq K
\text{ and $K$ is a maximal normal $\X$-subgroup}\}$;
in other words, $\nsup_G(H)$ is the \enquote{lowest level}
at which a maximal normal $\X$-subgroup of $G$ containing $H$ occurs.
\end{itemize}

Before we begin with the proof of Theorem~\ref{Thm:A} we record some preliminary  observations.

\begin{itemize}

\item[] {\bf Remark 1.} For $G$ and $H$ as above we have $\mathfrak{a}_G(H) = \mathfrak{a}_N(H)$,
where $N = \N_G(H)$. This is easily seen as every $\X$-subgroup of $G$
directly above $H$ (if there are any) normalises $H$.
\item[] {\bf Remark 2.}  It is clear that in~\ref{conj:a} we need to distinguish
between the cases $\mathfrak{a}_G(H) = 0$ 
and $\mathfrak{a}_G(H) \neq 0$ 
since if $H$ is
already a maximal $\X$-subgroup of $G$ 
then we have that $\mathfrak{a}_G(H) = 0$.
\item [] {\bf Remark 3.} If $H$ is an $\X$-subgroup of $G$ 
and has order $|H| = p^m$ 
then $\mathfrak{a}_G(H)$ is simply $\mathfrak{c}_{m + 1,G}(H)$ 
(which, in turn, will be $0$ 
in case $H$ is a maximal $\X$-subgroup of $G$).
\item [] {\bf Remark 4.} Assuming that $M$ is a maximal $\X$-subgroup of $G$
which contains $H$, it is clear that
$\sup_{G}(H) \leq |M| = \sup_{M}(H)$.
\end{itemize}

\vskip 1em
We are ready to prove Theorem~\ref{Thm:A}. 
The first step of the proof (\ref{conj:a} $\to$ \ref{conj:b}) 
resembles Yanovski's method of proof using   M\"{o}bius inversion
that can be found in~\cite{Yanovski}.

\begin{proofofA}
\ref{conj:a} $\to$ \ref{conj:b}
We induce on the index $|G:H|$.
In case $H=G$ we have (trivially) $\mathfrak{b}_G(H) = 1$ 
and thus $\mathfrak{b}_G(H) \equiv 1 \md{p}$. 
Observe that the exact same conclusion is reached if
$H$ is a maximal $\X$-subgroup of $G$, 
since $\mathfrak{b}_G(H) = 1$ 
in that case as well. 
We may therefore assume that $H$ is not a maximal $\X$-subgroup.
Denote by $\mathcal{M}$ the set of maximal $\X$-subgroups of $G$ so that
\[
\mathfrak{b}_G(H) = \left\lvert \{H \leq K \leq G: K \in \mathcal{M}\}\right\rvert.
\]
Now let $\mathcal{I}(G)$ be the set of subgroups of $G$ and
$\delta_G : \mathcal{I}(G) \to \{0,1\}$ be the indicator function
of the maximal $\X$-subgroups of $G$, i.e.
\[
\delta_{G}(K)=\left\{\begin{array}{l}
1, \text { if } K \in \mathcal{M} \\
0, \text { otherwise. }
\end{array}\right.
\]
Then $\mathfrak{b}_G(H) = \sum\limits_{H \leq K} \delta_{G}(K)$. 
Hence, by M\"{o}bius inversion, we get
\begin{equation}\label{Eq:First}
\delta_{G}(H) = \sum_{H \leq K} \mu(H, K) \mathfrak{b}_G(K).
\end{equation}
But $\mu(H, H) = 1$, while
\[
\mu(H, K) = \left\{\begin{array}{l}
(-1)^{r} p^{\binom{r}{2}}, \text { if } H \unlhd K \text { and } K / H \cong\left(C_{p}\right)^{r} \\
0, \text { otherwise }
\end{array}\right.
\]
(cf.~\cite[Prop. 2.4]{kt}).
Therefore equation~\eqref{Eq:First} becomes
\begin{equation}\label{Eq:Second}
\delta_{G}(H) = \mathfrak{b}_G(H) + \sum_{H < K} \mu(H, K) \mathfrak{b}_G(K) \equiv 
\mathfrak{b}_G(H) - \sum_{|K : H| = p} \mathfrak{b}_G(K) \md{p};
\end{equation}
the congruence above holds true due to the fact that 
$\mu(H,K) \equiv 0 \md{p}$ in case $|K:H| \geq p^2$. 
Equation~\eqref{Eq:Second} can now be rewritten as
\begin{equation}
\delta_{G}(H) \equiv \mathfrak{b}_G(H) - 
\sum_{|K : H| = p} \mathfrak{b}_G(K) \equiv
\mathfrak{b}_G(H) - \sum_{\substack{|K : H| = p \\ K \in \X}} 
\mathfrak{b}_G(K) \md{p},
\end{equation}
where the second congruence follows from 
the fact that $\mathfrak{b}_G(K)$ vanishes if $K \not\in \X$.
We conclude that 
\begin{equation}
\mathfrak{b}_G(H) \equiv
\sum_{\substack{|K : H| = p \\ K \in \X}} \mathfrak{b}_G(K) \md{p}
\end{equation}
since $H \not\in \mathcal{M}$ and thus $\delta_G(H) = 0$.
Notice that the number of terms in the sum above 
equals $\mathfrak{a}_G(H)$ and that for each $K \in \X$
such that $|K : H| = p$ we have $\mathfrak{b}_G(K) \equiv 1 \md{p}$
by the induction hypothesis applied to $K$.
Therefore 
\begin{equation}
\mathfrak{b}_{G}(H) \equiv 
\mathfrak{a}_G(H) \equiv 1 \md{p},
\end{equation}
where the last congruence holds by assumption.

\ref{conj:b} $\to$ \ref{conj:c} 
If $|H| = p^m$ the claim holds trivially for $\mathfrak{c}_{m,G}(H)$,
so we fix a $k$ with $p^m < p^k \leq \sup_G(H)$.
We let $\mathcal{M} = \{M_1, \ldots, M_s\}$ be the set of
maximal $\X$-subgroups of $G$ which contain $H$
and $\mathcal{K}=\{K_1,\ldots,K_t\}$ be the set of $\X$-subgroups
of $G$ which contain $H$ and have order $p^k$.
We are assuming that $s = \mathfrak{b}_G(H) \equiv 1 \md{p}$ 
and our goal is to prove that $t = \mathfrak{c}_{k,G}(H) \equiv 1 \md{p}$. 
Note that every subgroup in $\mathcal{M}$ contains some
subgroup in $\mathcal{K}$ and every subgroup in $\mathcal{K}$ is
contained in some subgroup in $\mathcal{M}$.

For each $i$ such that $1 \leq i \leq t$,
let $b_i$ be the number of subgroups in $\mathcal{M}$
which contain $K_i$. Observe that $b_i = \mathfrak{b}_G(K_i)$ 
and thus $b_i \equiv 1 \md{p}$ by assumption.
Similarly, for each $j \in \{1, \ldots, s\}$ let $c_j$
be the number of members of $\mathcal{K}$ which are
contained in $M_j$. 
Calculating the size of the set of ordered pairs $(i,j)$ 
such that $K_i \leq M_j$ in two ways yields
\begin{equation}\label{Eq:Third}
\sum_{i=1}^{t} b_i = \sum_{j=1}^{s} c_j.
\end{equation}

At this point, observe that $c_j = \mathfrak{c}_{k,M_j}(H)$.
Since $M_j$ is an $\X$-group, for every integer $k$ 
with $|H| \leq p^k \leq \sup_G(H) \leq |M|$
the quantity $\mathfrak{c}_{k,M_j}(H)$
is simply the number of subgroups of $M_j$ which have order
$p^k$ and contain $H$.
We now apply the standard result (cf.~\cite[Thm. 5.14]{bj1}) 
that the number of subgroups of a $p$-group at each level
containing a given subgroup is $\equiv 1 \md{p}$
to conclude that $c_j \equiv 1 \md{p}$.

We recall that $b_i \equiv 1 \md{p}$ for each $i \leq t$ 
and also $s \equiv 1 \md{p}$. 
Therefore the right-hand-side of~\eqref{Eq:Third} 
is congruent to $1 \md{p}$ and thus
\begin{equation}\label{Eq:Fourth}
1 \equiv \sum_{i=1}^{t} b_i \equiv t \md{p}.
\end{equation}
Since $k$ was arbitrary, we have established the claim.

\ref{conj:c} $\to$ \ref{conj:a} 
In case $H$ is a maximal
$\X$-subgroup of $G$, we have $\mathfrak{a}_G(H) = 0$ trivially 
and thus~\ref{conj:a} holds. 
Therefore we may assume that $\sup_G(H) > |H|$.
Let $|H| = p^m$. 
Then $\sup_G(H) \geq p^{m+1}$, so by assumption
\[
\mathfrak{a}_G(H) = \mathfrak{c}_{m + 1,G}(H) \equiv 1 \md{p}
\]
and thus~\ref{conj:a} is true in all cases. 
Our proof is complete.
\end{proofofA}

We continue with the proof of our
second main theorem.

\begin{proofofB}
Assume first that $G$ satisfies \ref{conj:a}, \ref{conj:b}, \ref{conj:c} in Theorem~\ref{Thm:A}. 

For the proof of \ref{Thm2a}, assume that 
$H$ is a maximal normal $\X$-subgroup $H$ of $G$. 
According to condition \ref{conj:a} 
either $\mathfrak{a}_G(H) = 0$ or $\mathfrak{a}_G(H) \equiv 1 \md{p}$. 
The latter is not possible however, since 
it is easy to see that if $\mathfrak{a}_G(H) \neq 0$
then $\mathfrak{a}_G(H)$ is congruent $\md{p}$ to the
number of normal $\X$-subgroups directly above it.
Thus if $H$ is a maximal normal $\X$-subgroup of $G$
it must be the case that $\mathfrak{a}_G(H) = 0$
and thus $H$ is a maximal $\X$-subgroup.

It remains to show \ref{conj:d}. 
Let $H \unlhd G$ with $|H|=p^m$. If $H$ is a maximal normal $\X$-subgroup then 
$\mathfrak{c}_{m,G}(H)= 1$ and  \ref{conj:d} holds trivially with $k=m$. 
So we may assume that $H$ is not a maximal normal $\X$-subgroup. Hence  $p^m < \nsup_G(H)$. Clearly for $k=m$ we have $\mathfrak{c}_{m,G}(H)=1$ and the claim holds. So we fix $k$ with $p^m < p^k \leq \nsup_G(H)$.

We look at the following sets of $\X$-subgroups of $G$ that contain $H$ 
\begin{enumerate}
 \item $\mathcal{M} = \{M_1, \ldots, M_t\}$ consists of 
maximal normal $\X$-subgroups of order $p^k$. So if $p^k < \nsup_G(H)$ we have $t=0$.
  \item $\mathcal{K} = \{K_1, \ldots, K_u\}$ consists of
non-maximal $\X$-subgroups of order $p^k$.
  \item $\mathcal{L}=\{L_1, \ldots, L_s\}$ consists of
non-normal maximal $\X$-subgroups of order $ \leq p^k$. The first $s_1$ of them $\{ L_1, \cdots , L_{s_1}\}$ 
are of order exactly $p^k$.
 \item $\mathcal{T} = \{T_1, \ldots, T_w\}$ consists of 
maximal $\X$-subgroups of order $> p^k$. 
\end{enumerate}

First observe that if $S$ is an $\X$-subgroup of $G$ that contains $H$ then the same holds for any $G$-conjugate $S^g$ of $S$. In addition $S$ is maximal if and only if $S^g$ is maximal. Hence the number $s_1$ of non-normal maximal $\X$-subgroups of order $p^k$ is a multiple of $p$ as it is a union of conjugacy classes of subgroups that are not normal. The same holds for $s$ the total number of non-normal maximal $\X$-subgroups of order $ \leq p^k$ that contain $H$.

Furthermore, the definition of $\mathfrak{c}_{k,G}(H)$ implies that $\mathfrak{c}_{k,G}(H)= u+t+s_1$. Hence the previous observation implies 
\begin{equation}\label{Eq:T2.1}
 \mathfrak{c}_{k,G}(H)= u+t+s_1 \equiv u+t \md{p}
\end{equation}

Now the total number $b_G(H)$ of maximal $\X$-subgroups that contain $H$ is $s+t+w$. In addition, $b_G(H)$ is congruent to 1 $\md{p}$, by hypothesis. Hence 
\begin{equation}\label{Eq:T2.2}
 1\equiv \mathfrak{b}_{G}(H)= s+t+w \equiv t+w \md{p}
\end{equation}
where the last equivalence holds because $s \equiv 0 \md{p}$.

Let $b_i=| \{ T_j \in \mathcal{T} \mid K_i \subseteq T_j\} |$
and $c_j= | \{ K_l \in \mathcal{K} \mid K_l \subseteq T_j\} |$ for all $i =1, \cdots , u$ 
and all $j=1, \cdots , w$. The fact that every maximal normal $\X$-subgroup is a maximal subgroup implies that 
none of the $T_j$ contains any of the $M_i$. So $b_i, c_j \neq 0$ and in addition, 
\begin{equation}\label{Eq:T2.3}
\sum_{i=1}^{u} b_i = \sum_{j=1}^{w} c_j.
\end{equation}
Observe that $b_i= \mathfrak{b}_G(K_i)$ and thus $b_i \equiv 1 \md{p}$ by assumption, since $G$ satisfies Condition B. Also,  $c_j = \mathfrak{c}_{k,T_j}(H)$. But $T_j$ is an $\X$-group, thus the number of $\X$-subgroups of $T_j$ of order $p^k$ containing $H$ is equal to the total number of 
subgroups of $T_j$ of order $p^k$ containing $H$. The latter is $\equiv 1 \md{p}$ according to Theorem 5.14 in ~\cite[Thm. 5.14]{bj1}).
Hence 
\[
c_j = \mathfrak{c}_{k,T_j}(H) \equiv 1 \md{p}.
\]
We conclude that equation \eqref{Eq:T2.3} implies 
\begin{equation*}
u \equiv w \md{p}
\end{equation*}
This along with equations \eqref{Eq:T2.1} and \eqref{Eq:T2.2} provide 
\[
\mathfrak{c}_{k,G}(H) \equiv u+t \equiv w+t \equiv \mathfrak{b}_{G}(H) \equiv 1 \md{p}
\] 

For the other direction, assume that $G$ satisfies \ref{Thm2a} and \ref{conj:d}. 
We will show, using induction on $|G:H|$, that \ref{conj:a} holds. 
The base case $G = H $ is 
trivially true. If $\N_G(H) < G$, then since
$\mathfrak{a}_G(H) = \mathfrak{a}_N(H)$,
where $N = \N_G(H)$, the claim follows by the induction
hypothesis applied to $H$ with respect to $N$.
We may thus assume that $H$ is normal in $G$.
Suppose $|H| = p^k$. If $\nsup_G(H) > |H|$, 
then $\mathfrak{c}_{k+1,G} \equiv 1 \md{p}$
whence $\mathfrak{a}_G(H) \equiv 1 \md{p}$.
There only remains to consider the case
$\nsup_G(H) = |H|$, i.e. when $H$ is a maximal
normal $\X$-subgroup of $G$. But then the
assumption that every maximal normal $\X$-subgroup of $G$ 
is a maximal $\X$-subgroup of $G$ 
forces $H$ to be a maximal $\X$-subgroup of $G$ 
and thus $\mathfrak{a}_G(H) = 0$.
This completes the induction and proves the claim.
\end{proofofB}

\section{Abelian Subgroups}
The existence of (maximal) normal abelian subgroups of fixed order 
in a finite $p$-group $G$ has been studied by several authors, 
see for example \cite{gl1}, \cite{gl2}, \cite{berk}. 
In this section we wish to explore how our Theorems~\ref{Thm:A} and~\ref{Thm:B}
can be applied to the class of abelian groups 
of given order and bounded exponent. 

We begin with the following.

\begin{theoremno}[{\cite[Thm. 1]{berk}}]
Let $A < B \leq G$, where $A, B$ are abelian subgroups 
of a $p$-group $G$, $|B : A| = p$, $\exp(B) \leq p^{k}$ 
and $p^{k} > 2$. Let $\mathcal{A}$ be the set of all 
abelian subgroups $T$ of $G$ such that $A < T$, $|T : A| = p$ 
and $\exp (T) \leq p^{k}$. 
Then $|\mathcal{A}| \equiv 1 \md{p}$.
\end{theoremno}

This theorem of Berkovich says, essentially, that if $p$, $k$ 
are such that $p^k > 2$, $G$ is a finite $p$-group 
and $H$ is a subgroup of $G$, then either 
$\mathfrak{a}_G(H) = 0$ or $\mathfrak{a}_G(H) \equiv 1 \md{p}$
relative to the class $\X$ of finite abelian $p$-groups
of exponent $\leq p^k$.

Since~\ref{conj:a} in Theorem~\ref{Thm:A} holds for the
class $\X$ of finite abelian $p$-groups of exponent $\leq p^k$
for some fixed choice of $k$ (subject only to the condition 
$p^k > 2$), we see that~\ref{conj:b} holds as well and thus
 we recover the main result in the  recent paper of Isaacs and Yanovski. 

\begin{prevcorollary}[{\cite[Thm. C]{isaacs}}]\label{Cor:A}
Let $P$ be a $p$-group, and suppose that $e>2$ is a power of $p$. 
Also, let $H \leq P$ be an abelian subgroup with exponent dividing $e$, 
and let $n$ be the number of subgroups $A$ of $P$ that contain $H$ 
and that are maximal with respect to the property 
that $A$ is abelian and has exponent dividing $e$. 
Then $n \equiv 1 \md{p}$.
\end{prevcorollary}

Moreover, given a finite $p$-group $G$ of exponent $p^k$
we let $\X$ be the class of finite abelian $p$-groups
of exponent $\leq p^k$. 
Then~\ref{conj:a} in Theorem~\ref{Thm:A} holds for the
class $\X$ and thus~\ref{conj:b} holds.
The maximal $\X$-subgroups which contain a given $\X$-subgroup
$H$ of $G$ are then simply the maximal abelian subgroups
of $G$ containing $H$. The cardinality of the first set
is our quantity $\mathfrak{b}_G(H) \equiv 1 \md{p}$
and thus we recover Yanovski's result mentioned 
in the introduction.

A celebrated theorem due to Burnside asserts that 
if $G$ is a $p$-group of order $p^n$ and $A$ is a maximal
abelian normal subgroup of $G$, where $|A| = p^s$, then
$n \leq \frac{s(s+1)}{2}$. 
Since this holds for all maximal abelian normal subgroups of $G$, 
we see that $\nsup_G(1) \geq p^s$ with $s$
the least positive integer for which $n \leq \frac{s(s+1)}{2}$ holds
relative to the class of finite abelian groups $\X = \mathfrak{A}$.

Now every finite $p$-group satisfies Theorem~\ref{Thm:A} 
with $\X = \mathfrak{A}$ and thus every finite $p$-group $G$ 
satisfies Theorem~\ref{Thm:B}.
If $k \leq s$ therefore, then $k$ lies in the range 
for which \ref{conj:d} of Theorem~\ref{Thm:B} works.
But 
\[
\frac{k(k-1)}{2} \leq \frac{s(s-1)}{2} < n
\]
where the last inequality is a consequence of the minimality of $s$.
We deduce that $\binom{k}{2} < n$ and thus we recover the following classic result due to Miller.

\begin{prevcorollary}[{\cite[Thm. 13.12]{bj1}}]
Suppose that $G$ is a group of order $p^{n}$, 
$n, k \in \mathbb{N}$ and $n > \binom{k}{2}$. 
Then the number of abelian subgroups of order $p^{k}$ in $G$ 
is congruent to $1 \md{p}$.
\end{prevcorollary}

The next theorem provides an  analogue of Burnside's theorem 
for maximal abelian normal subgroups of given exponent.

In order to prove the theorem, we need the following auxiliary result 
which is essentially Lemma 1 in~\cite{laffey80}. 
We follow the notation in~\cite{laffey80} and write 
\[
\Omega (G) = \begin{cases}
\Omega_2(G), \text{\quad if } p=2 \\
\Omega_1(G), \text{\quad if } p>2
\end{cases}
\]
for any $p$-group $G$. 
We also write $d(G)$ for the minimum number of generators of $G$.

\begin{prevlemma}\label{Lem:laffey}
Let $C$ be a $p$-group of exponent $\leq p^e$. 
Assume further that $\Omega_s (C) \leq \zent(C)$. 
Then 
$|C|/|\Omega_s(C)|\leq |\Omega_1(C)|^{e-s}$ 
for every integer $s$ with $p^s >2$.
\end{prevlemma}

\begin{proofof}
First note that $\Omega(C) \leq \Omega_s(C) \leq Z(C)$ 
for every integer $s$ with $p^s >2$.
Hence Corollary 1 of Lemma 1 in~\cite{laffey74} implies that 
$\Omega_i(C) = \{x \in C \mid x^{p^i} = 1\}$ for every $i \geq 0$. 
Therefore $\Omega_{i+1}(C)/\Omega_i(C)$ 
is an elementary abelian group for all $i \geq 1$. 
So
\[
\left\lvert\Omega_{i+1}(C)/\Omega_i(C)\right\rvert \leq 
|\Omega_{i+1}(C)/\fratt(\Omega_{i+1}(C))| \leq 
d(\Omega_{i+1}(C)).
\]
In addition, for every $i \geq 1$, we have 
\[
\Omega( \Omega_i(C)) \leq \Omega (C) \leq \Omega_s(C) \leq \zent(C).
\]
Applying Corollary 2 in \cite{laffey74} on $\Omega_i(C)$ we get 
\[
d(\Omega_{i}(C)) \leq d(\Omega_1(\Omega_i(C)) = d(\Omega_1(C)),
\]
where the last equality follows from the fact that 
$\Omega_1(\Omega_i(C)) = \Omega_1(C)$. 

Now notice that if $d(\Omega_1(C)) = t$ then $|\Omega_1(C)|=p^t$. 
Hence for every $i \geq s$ we get $|\Omega_{i+1}(C)/\Omega_i(C)| \leq p^t$.
As $C = \Omega_e(C)$ we need to repeat the argument $(e - s)$ times to get 
$|C|/|\Omega_s(C)| \leq p^{t(e-s)}$, and the result follows.
\end{proofof}

We are now ready to prove our last
main theorem.

\begin{proofofC}
We denote by $d = d(A)$ the rank of the subgroup $A$ and begin with
the observation that $d \leq r - s + 1$. 
This is so because 
$A = B \times D$ for some subgroups $B$, $D$ of $A$, 
where $B$ is cyclic of order $p^s$. 
Thus
\[
d(A) = d(D) + 1 \leq r - s + 1,
\]
the inequality being a consequence of $p^{d(D)} \leq |D|$.

Our strategy will be to write 
\[
|G| = |G:C||C:A||A|,
\]
where $C \coloneqq \cent_G(A)$,
and to obtain bounds for the indices involved.
As regards the quantity $|G:C|$, the $N/C$ Theorem
tells us that $|G:C| \leq \left\lvert \aut(A)\right\rvert_p$, 
the $p$-part of $\left\lvert \aut(A)\right\rvert$.
But 
\[
\left\lvert \aut(A)\right\rvert_p \leq 
p^{\frac{1}{2}d(2r - d - 1)} \leq 
p^{\frac{1}{2}(r - s + 1)[2r - (r - s + 1) - 1]} = 
p^{\frac{1}{2}(r - s + 1)(r + s - 2)}.
\]
The first inequality follows from the fact that 
$\left\lvert \aut(A)\right\rvert$ divides 
$p^{d(r-d)} \left\lvert \gl_d(p) \right\rvert$ 
(cf.~\cite[Satz 3.19]{endliche}).
To see why the second inequality is true, 
let $f(x) = \frac{1}{2}x(2r - x - 1)$ and note that
\[
f(d) - f(r - s + 1) = \frac{1}{2} 
\overbrace{(2 + d - r - s)}^{\leq 0} 
\overbrace{(r - s + 1 - d)}^{\geq 0} \leq 0\, ;
\]
we have $2 + d - r - s \leq 0$ since $2 \leq s$ and
$d \leq r$, while $r - s + 1 - d \geq 0$ by our initial observation.

Regarding the quantity $|C:A|$, we use Lemma~\ref{Lem:laffey}
to bound it. 
By a theorem of Alperin~\cite{alperin} $\Omega_s(C) = A$
and thus $|C : A| \leq |\Omega_1(C)|^{e-s}$ 
by~Lemma~\ref{Lem:laffey}.
But 
\[
|\Omega_1(C)| = p^{d(\Omega_1(C))} 
  \leq p^{d(A)} 
  \leq p^{r - s + 1}.
\]
The equality is a consequence of the fact that
\[
\Omega_1(C) \leq \Omega_s(C) = A \leq \zent(C).
\]
The first inequality follows from the containment
$\Omega_1(C) \leq \Omega_s(C) = A$ and the fact that $A$
is abelian.
Thus
\[
|C : A| \leq \left(p^{r - s + 1}\right)^{e-s} = 
p^{(e - s)(r - s + 1)}.
\]
Then 
\begin{align}
|G| = |G:C||C:A||A| &\leq p^{\frac{1}{2}(r - s + 1)(r + s - 2)}
    p^{(e - s)(r - s + 1)}
    p^r \\
   &= p^{\frac{1}{2}(r - s + 1)(r - s + 2e - 2) +r},
\end{align}
which is clearly equivalent to what we wanted to prove.

Now,  the inequality 
\[
2n \leq r^2 +r(2e-2s+1) +(s-1)(s-2e+2).
\]
along with the fact that for every $e \geq 2$ we get $(s-1)(s-2e+2)\leq 0$ 
(since $s \leq e \leq 2e-2$) implies 
\[
2n \leq r^2 +r(2e-2s+1).
\]
Solving the quadratic for $r$ yields 
$r \geq -(e-s +1/2) + \sqrt{(e-s+1/2)^2 +2n}$. This completes the proof of the theorem.
\end{proofofC}

In Theorem~\ref{Thm:C} we have not addressed the case $s = 1$
which corresponds to the class of elementary abelian subgroups of a $p$-group.
But this case has already been handled by Laffey in~\cite{laffey80}. 
The main result of that paper is the following.

\begin{theoremno}[{\cite[Thm. 1]{laffey80}}]
Let $G$ be a finite $p$-group, 
and suppose that $|G| = p^n$, $\exp(G) = p^e$ and $p^r = \min\{|A|\}$,
where the minimum is taken over 
all maximal elementary abelian normal subgroups $A$ of $G$.
Then \[n \leq \begin{cases}
\begin{aligned}
&\frac{1}{2}r(3r - 1) + er, &\text{ if } p =2,\\
&\frac{1}{2}r(r - 1) + er, &\text{ if } p > 2.
\end{aligned}
\end{cases}
\]
\end{theoremno}

The quantity $p^r$ in Laffey's result is precisely our $\nsup_G(1)$
relative to the class $\mathfrak{Y}$ of elementary abelian $p$-groups.
Moreover, Laffey's inequality implies that (relative to the class $\mathfrak{Y}$)
\[
\nsup_G(1) = r \geq \begin{cases}
\begin{aligned}
-y +&\sqrt{ y^2 + \frac{2n}{3}} &\text{ if } p=2, \\ 
 -y + &\sqrt{y^2 + 2n} &\text{ if } p>2
\end{aligned}
\end{cases}
\]
where 
\[
y= \begin{cases}
\frac{1}{6}(2e-1) \text{ if } p = 2\\
\frac{1}{2}(2e-1) \text{ if } p > 2 .
\end{cases}
\]

Thus we can apply Theorem~\ref{Thm:B}~\ref{conj:d} to the class $\mathfrak{Y}$ (of elementary abelian subgroups) when $H = 1$ to deduce the following.

\begin{prevcorollary}
Suppose that $G$ is a group of order $p^{n}$ and exponent $p^e$, while $y$ is the function of $e$ defined above. 
Then for all $1\leq k$ such that 
\[
 k \leq \begin{cases}
 \begin{aligned}
   -y + &\sqrt{y^2 +\frac{2n}{3}} &\text{ if } p =2, \\
   -y + &\sqrt{y^2 + 2n} &\text{ if } p >2
 \end{aligned}
  \end{cases}
\]
the number of elementary abelian subgroups of order $p^k$ is\, $\equiv 1 \md{p}$. In particular, there exists at least one
elementary abelian normal subgroup of order $p^k$ for all such $k$.
\end{prevcorollary}

Additionally, solving the inequality
\[
-y + \sqrt{y^2 + 2n} \geq m
\]
allows one to conclude that a finite $p$-group $G$ of odd order
$p^n$ and exponent $p^e$ is guaranteed to have an elementary abelian normal subgroup
of order $p^k$ for all $k \leq m$ provided that
\[
e \leq \frac{n}{m} -\frac{m-1}{2}
\]
holds.
Compare the corollary above with Theorems 6 and 7 in~\cite{berk}.

\section{Some problems and questions}
In this final section we outline some open problems 
and questions, hoping thereby to stimulate further research on this topic.

The only instances of classes $\X$ for which Theorem~\ref{Thm:A}
holds universally (i.e. for all finite $p$-groups) are subclasses
of $\mathfrak{A}$, the class of finite abelian groups.
One natural problem here is the following: 

\begin{problem}
Find other classes $\X$
which are not subclasses of $\mathfrak{A}$ such that all finite
$p$-groups satisfy Theorem~\ref{Thm:A} with respect to $\X$ or prove that no such classes exist.
\end{problem}

\begin{question}
What can we conclude if a group $G$ 
satisfies $\sup_G(1) = \nsup_G(1)$?
\end{question}

Finally, is a certain \enquote{converse} true?
That is: 
\begin{question}
Suppose that $\X$ enjoys the following
property: for all $p$-groups $G$, if a subgroup is 
a maximal normal $\X$-subgroup of $G$ then it is also
a maximal $\X$-subgroup of $G$. 
Is it then true that every $p$-group $G$ satisfies
the conclusions of Theorem~\ref{Thm:A} with respect to $\X$?
\end{question}

\section*{Acknowledgement}
Both  authors  would like to thank Prof. George Glauberman  for his valuable remarks.

\bibliographystyle{amsalpha}
\bibliography{Bibliography}
\end{document}